\documentclass[12pt]{amsart}
\usepackage{amsmath,amsthm,amsfonts,amssymb,mathrsfs}
\date{\today}

\usepackage[cp1251]{inputenc}         % Кодова сторінка Windows1251

\usepackage[ukrainian]{babel} % Підтримка кирилиці
\usepackage{amsmath,amsthm,amsfonts,amssymb,mathrsfs}
\usepackage{eucal}

\input xymatrix
\xyoption{all}

\usepackage{color}

\usepackage{amssymb,cite}
\usepackage[colorlinks,plainpages,citecolor=magenta, linkcolor=blue, backref]{hyperref}%,bookmarksnumbered,

\usepackage{hyperref}

\usepackage{tikz}

\newtheorem{theorem}{Теорема}%[section]
\newtheorem{question}{Питання}
\newtheorem{proposition}{Твердження}

\newtheorem{lemma}{Лема}
\theoremstyle{definition}
\newtheorem{example}{Приклад}%[section]
\newtheorem{remark}{Зауваження}%[section]

\newtheorem{definition}[theorem]{Означення}%[section]

\begin{document}

\title[Про автоморфiзми напiвгрупи $\boldsymbol{B}_{\omega}^{\mathscr{F}}$  ...]{Про автоморфiзми напiвгрупи $\boldsymbol{B}_{\omega}^{\mathscr{F}}$  у випадку сiм'ї $\mathscr{F}$ непорожнiх iндуктивних пiдмножин у  ${\omega}$}

\author[Олег~Гутік, Микола Михаленич]{Олег~Гутік, Микола Михаленич}
\address{Механіко-математичний факультет, Львівський національний університет ім. Івана Франка, Університецька 1, Львів, 79000, Україна}
\email{oleg.gutik@lnu.edu.ua,
ogutik@gmail.com, myhalenychmc@gmail.com}



\begin{abstract}
Нехай $\mathscr{F}$ --- сім'я непорожнiх iндуктивних пiдмножин у  ${\omega}$.
Доведено, що iн'\-cктивний ендоморфізм $\varepsilon$ напiвгрупи $\boldsymbol{B}_{\omega}^{\mathscr{F}}$ є тотожним відображенням тоді і тільки тоді, коли $\varepsilon$ має три різні нерухомі точки, що еквівалентно існуванню неідемпотентного елемента $(i,j,[p))\in\boldsymbol{B}_{\omega}^{\mathscr{F}}$ такого, що $(i,j,[p))\varepsilon=(i,j,[p))$.

\bigskip
\noindent
\emph{Oleg Gutik, Mykola Mykhalenych, \textbf{On automorphisms of the semigroup $\boldsymbol{B}_{\omega}^{\mathscr{F}}$ in the case when the family $\mathscr{F}$ consists of nonempty inductive subsets of ${\omega}$}.}

\smallskip
\noindent
Let $\mathscr{F}$ be a family of nonempty inductive subsets of ${\omega}$.
It is proved that an injective endomorphism $\varepsilon$ of the semigroup $\boldsymbol{B}_{\omega}^{\mathscr{F}}$ is the transformation if and only if  $\varepsilon$ has three distinct fixed points, which is equivalent to existence non-idempotent element $(i,j,[p))\in\boldsymbol{B}_{\omega}^{\mathscr{F}}$ such that  $(i,j,[p))\varepsilon=(i,j,[p))$.
\end{abstract}

\maketitle

%\section{Термінологія та означення}

%\section{Вступ}\label{section-1}

У цій праці ми користуємося термінологією з монографій \cite{Clifford-Preston-1961, Clifford-Preston-1967, Lawson-1998, Petrich-1984}.
Надалі у тексті множину невід'ємних цілих чисел  позначатимемо через $\omega$. Для довільного числа $k\in\omega$ позначимо
${[k)=\{i\in\omega\colon i\geqslant k\}.}$

Нехай $\mathscr{P}(\omega)$~--- сім'я усіх підмножин у $\omega$.
Для довільних $F\in\mathscr{P}(\omega)$ i цілого числа  $n$ приймемо
\begin{equation*}
n+F=\{n+k\colon k\in F\}, \qquad \hbox{якщо} \; F\neq\varnothing
\end{equation*}
i $n+F=\varnothing$, якщо $F=\varnothing$.
Будемо говорити, що непорожня підсім'я  $\mathscr{F}\subseteq\mathscr{P}(\omega)$ є \emph{${\omega}$-замкненою}, якщо $F_1\cap(-n+F_2)\in\mathscr{F}$ для довільних $n\in\omega$ та $F_1,F_2\in\mathscr{F}$.

Підмножина $A$ в $\omega$ називається \emph{індуктивною}, якщо з  $i\in A$ випливає, що $i+1\in A$. Очевидно, що $\varnothing$ --- індуктивна множина в $\omega$.

\begin{remark}\label{remark-2.2}
\begin{enumerate}
  \item\label{remark-2.2(1)} За лемою 6 з \cite{Gutik-Mykhalenych-2020} непорожня множина $F\subseteq \omega$ є індуктивною в $\omega$ тоді і лише тоді, коли $(-1+F)\cap F=F$.
  \item\label{remark-2.2(2)} Оскільки множина $\omega$ зі звичайним порядком є цілком впорядкованою, то для кожної непорожньої індуктивної множини $F$ у $\omega$ існує невід'ємне ціле число $n_F\in\omega$ таке, що $[n_F)=F$.
  \item\label{remark-2.2(3)} З \eqref{remark-2.2(2)} випливає, що перетин довільної скінченної кількості непорожніх індуктивних підмножин у $\omega$ є непорожньою індуктивною підмножиною в $\omega$.
\end{enumerate}
\end{remark}

Надалі, якщо $(X,\le)$~--- частково впорядкована множина та $x\in X$, то позначимо
\begin{equation*}
  \uparrow_{\le}{x}=\{y\in X\colon x\le y\}.
\end{equation*}

Якщо $S$~--- напівгрупа, то її підмножина ідемпотентів позначається через $E(S)$.  На\-пів\-гру\-па $S$ називається \emph{інверсною}, якщо для довільного її елемента $x$ існує єдиний елемент $x^{-1}\in S$ такий, що $xx^{-1}x=x$ та $x^{-1}xx^{-1}=x^{-1}$ \cite{Petrich-1984, Vagner-1952}. В інверсній напівгрупі $S$ вище означений елемент $x^{-1}$ називається \emph{інверсним до} $x$. \emph{В'язка}~--- це напівгрупа ідемпотентів, а \emph{напівґратка}~--- це комутативна в'язка. %Надалі через $(\mathscr{P}_{\!\infty}(X),\cup)$ по\-зна\-ча\-ти\-ме\-мо \emph{вільну напівґратку} з одиницею над не\-порож\-ньою множиною $X$, тобто множину усіх скінченних (включно з по\-рож\-ньою) підмножин множини $X$ з операцією ``об'єднання''.

%Якщо $S$ --- напівгрупа, то ми позначатимемо відношення Ґріна на $S$ через $\mathscr{R}$, $\mathscr{L}$, $\mathscr{D}$, $\mathscr{H}$ і $\mathscr{J}$ (див. означення в \cite[\S2.1]{Clifford-Preston-1961} або \cite{Green-1951}). Напівгрупа $S$ називається \emph{простою}, якщо $S$ не містить власних двобічних ідеалів, тобто $S$ складається з одного $\mathscr{J}$-класу, і \emph{біпростою}, якщо $S$ складається з одного $\mathscr{D}$-класу.

%Відношення еквівалентності $\mathfrak{K}$ на напівгрупі $S$ називається \emph{конгруенцією}, якщо для елементів $a$ та $b$ напівгрупи $S$ з того, що виконується умова $(a,b)\in\mathfrak{K}$ випливає, що $(ca,cb), (ad,bd) \in\mathfrak{K}$, для довільних $c,d\in S$. Відношення $(a,b)\in\mathfrak{K}$ ми також будемо записувати $a\mathfrak{K}b$, і в цьому випадку будемо говорити, що \emph{елементи $a$ i $b$ є $\mathfrak{K}$-еквівалентними}.

%Конгруенція $\mathfrak{K}$ на напівгрупі $S$ називається \emph{груповою}, якщо фактор-напівгрупа $S/\mathfrak{K}$ ізоморфна деякій групі $G$ \cite{Clifford-Preston-1961}. Нагадаємо \cite{Lawson-1998, Petrich-1984}, що на кожній ін\-версній напівгрупі $S$ існує \emph{найменша} (\emph{мінімальна}) \emph{групова конгруенція} $\boldsymbol{\sigma}$ і вона визначається так:
%\begin{equation*}
%  s\boldsymbol{\sigma}t \qquad \Longleftrightarrow \qquad es=et \quad \hbox{для деякого} \quad e\in E(S).
%\end{equation*}

Якщо $S$~--- напівгрупа, то на $E(S)$ визначено частковий порядок:
$
e\preccurlyeq f
$   тоді і лише тоді, коли
$ef=fe=e$.
Так означений частковий порядок на $E(S)$ називається \emph{при\-род\-ним}.

Означимо відношення $\preccurlyeq$ на інверсній напівгрупі $S$ так:
$
    s\preccurlyeq t
$
тоді і лише тоді, коли $s=te$, для деякого ідемпотента $e\in S$. Так означений частковий порядок назива\-єть\-ся \emph{при\-род\-ним част\-ковим порядком} на інверсній напівгрупі $S$~\cite{Vagner-1952}. Очевидно, що звуження природного часткового порядку $\preccurlyeq$ на інверсній напівгрупі $S$ на її в'язку $E(S)$ є при\-род\-ним частковим порядком на $E(S)$.

Нагадаємо (див.  \cite[\S1.12]{Clifford-Preston-1961}, що \emph{біциклічною напівгрупою} (або \emph{біциклічним моноїдом}) ${\mathscr{C}}(p,q)$ називається напівгрупа з одиницею, породжена двоелементною мно\-жи\-ною $\{p,q\}$ і визначена одним  співвідношенням $pq=1$.
%Біциклічна на\-пів\-група відіграє важливу роль у теорії на\-півгруп. Так, зокрема, класична теорема О.~Ан\-дерсена \cite{Andersen-1952}  стверджує, що {($0$-)}прос\-та напівгрупа з (ненульовим) ідем\-по\-тен\-том є цілком {($0$-)}прос\-тою тоді і лише тоді, коли вона не містить ізоморфну копію бі\-циклічного моноїда. Різні розширення та узагальнення біциклічного моноїда вводилися раніше різ\-ни\-ми авторами \cite{Fortunatov-1976, Fotedar-1974, Fotedar-1978, Gutik-Pagon-Pavlyk=2011, Warne-1967}. Такими, зокрема, є конструкції Брука та Брука--Рейлі занурення напівгруп у прості та описання інверсних біпростих і $0$-біпростих $\omega$-напівгруп \cite{Bruck-1958, Reilly-1966, Warne-1966, Gutik-2018}.

\begin{remark}\label{remark-10}
Легко бачити, що біциклічний моноїд ${\mathscr{C}}(p,q)$ ізоморфний напівгрупі, заданій на множині $\boldsymbol{B}_{\omega}=\omega\times\omega$ з напівгруповою операцією
\begin{equation*}
  (i_1,j_1)\cdot(i_2,j_2)=(i_1+i_2-\min\{j_1,i_2\},j_1+j_2-\min\{j_1,i_2\})=
\left\{
  \begin{array}{ll}
    (i_1-j_1+i_2,j_2), & \hbox{якщо~} j_1\leqslant i_2;\\
    (i_1,j_1-i_2+j_2), & \hbox{якщо~} j_1\geqslant i_2.
  \end{array}
\right.
\end{equation*}
Надалі під біциклічною напівгрупою, чи біциклічним моноїдом, ми будемо розуміти напівгрупу $\boldsymbol{B}_{\omega}$.
\end{remark}

У праці \cite{Gutik-Mykhalenych-2020} введено алгебраїчні розширення $\boldsymbol{B}_{\omega}^{\mathscr{F}}$ біциклічного моноїда для довільної $\omega$-замк\-не\-ної сім'ї $\mathscr{F}$ підмножин в $\omega$, які узагальнюють біциклічний моноїд, зліченну напівгрупу матричних одиниць і деякі інші комбінаторні інверсні напівгрупи.

Нагадаємо цю конструкцію.
Нехай $\boldsymbol{B}_{\omega}$~--- біциклічний моноїд і  $\mathscr{F}$ --- непорожня ${\omega}$-замкнена підсім'я в  $\mathscr{P}(\omega)$. На множині $\boldsymbol{B}_{\omega}\times\mathscr{F}$ озна\-чимо бінарну операцію ``$\cdot$''  формулою
\begin{equation}\label{eq-1.1}
  (i_1,j_1,F_1)\cdot(i_2,j_2,F_2)=
  \left\{
    \begin{array}{ll}
      (i_1-j_1+i_2,j_2,(j_1-i_2+F_1)\cap F_2), & \hbox{якщо~} j_1<i_2;\\
      (i_1,j_2,F_1\cap F_2),                   & \hbox{якщо~} j_1=i_2;\\
      (i_1,j_1-i_2+j_2,F_1\cap (i_2-j_1+F_2)), & \hbox{якщо~} j_1>i_2.
    \end{array}
  \right.
\end{equation}

У \cite{Gutik-Mykhalenych-2020} доведено, якщо сім'я  $\mathscr{F}\subseteq\mathscr{P}(\omega)$ є ${\omega}$-замкненою, то $(\boldsymbol{B}_{\omega}\times\mathscr{F},\cdot)$ є напівгрупою. Ми в \cite{Gutik-Mykhalenych-2020} доводимо, що $\boldsymbol{B}_{\omega}^{\mathscr{F}}$ є комбінаторною інверсною напівгрупою, а також описано відношення Ґріна, частковий природний порядок на напівгрупі $\boldsymbol{B}_{\omega}^{\mathscr{F}}$ та її множину ідемпотентів. Також, у \cite{Gutik-Mykhalenych-2020} доведено критерії  простоти, $0$-простоти, біпростоти та $0$-біпростоти напівгрупи $\boldsymbol{B}_{\omega}^{\mathscr{F}}$, і вказано умови, коли $\boldsymbol{B}_{\omega}^{\mathscr{F}}$ містить одиницю, ізоморфна біциклічному моноїду або зліченній напівгрупі матричних одиниць.

Зауважимо, що у \cite{Gutik-Pozdnyakova-2021??} отримано подібні результати до \cite{Gutik-Mykhalenych-2020} у випадку розширення $\boldsymbol{B}_{\mathbb{Z}}^{\mathscr{F}}$ розширеної біциклічної напівгрупи $\boldsymbol{B}_{\mathbb{Z}}$ для довільної $\omega$-замк\-не\-ної сім'ї $\mathscr{F}$ підмножин в $\omega$, а в \cite{Gutik-Pozdniakova=2023} описано група автоморфізмів напівгрупи $\boldsymbol{B}_{\mathbb{Z}}^{\mathscr{F}}$.

Припустимо, що ${\omega}$-замкнена сім'я $\mathscr{F}\subseteq\mathscr{P}(\omega)$ містить порожню множину $\varnothing$, то з означення напівгрупової операції $(\boldsymbol{B}_{\omega}\times\mathscr{F},\cdot)$ випливає, що множина
$ %\begin{equation*}
  \boldsymbol{I}=\{(i,j,\varnothing)\colon i,j\in\omega\}
$ %\end{equation*}
є ідеалом напівгрупи $(\boldsymbol{B}_{\omega}\times\mathscr{F},\cdot)$.

\begin{definition}[\!\!{\cite{Gutik-Mykhalenych-2020}}]\label{definition-1.1}
Для довільної ${\omega}$-замкненої сім'ї $\mathscr{F}\subseteq\mathscr{P}(\omega)$ означимо
\begin{equation*}
  \boldsymbol{B}_{\omega}^{\mathscr{F}}=
\left\{
  \begin{array}{ll}
    (\boldsymbol{B}_{\omega}\times\mathscr{F},\cdot)/\boldsymbol{I}, & \hbox{якщо~} \varnothing\in\mathscr{F};\\
    (\boldsymbol{B}_{\omega}\times\mathscr{F},\cdot), & \hbox{якщо~} \varnothing\notin\mathscr{F}.
  \end{array}
\right.
\end{equation*}
\end{definition}

У \cite{Gutik-Lysetska=2021, Lysetska=2020} досліджено алгебраїчну структуру напівгрупи $\boldsymbol{B}_{\omega}^{\mathscr{F}}$ у випадку, коли ${\omega}$-замкнена сім'я $\mathscr{F}$ складається з атомарних підмножин (одноточкових підмножин і порожньої множини) в ${\omega}$. Зокрема доведено, що за виконання таких умов на сім'ю $\mathscr{F}$ напівгрупа  $\boldsymbol{B}_{\omega}^{\mathscr{F}}$ ізоморфна піднапівгрупі ${\omega}$-розширення Брандта напівгратки $(\omega,\min)$. Також, у \cite{Gutik-Lysetska=2021, Lysetska=2020} досліджувались топологізація напівгрупи $\boldsymbol{B}_{\omega}^{\mathscr{F}}$, близькі до компактних трансляційно неперервні топології на $\boldsymbol{B}_{\omega}^{\mathscr{F}}$ та замикання напівгрупи $\boldsymbol{B}_{\omega}^{\mathscr{F}}$ у напівтопологічних напівгрупах.

Із зауваження~\ref{remark-2.2}\eqref{remark-2.2(3)} випливає, якщо сім'я $\mathscr{F}_0$ склада\-єть\-ся з індуктивних в $\omega$ підмножин і міс\-тить порожню множину $\varnothing$ як елемент, то для сім'ї $\mathscr{F}=\mathscr{F}_0\setminus\{\varnothing\}$ множина  $\boldsymbol{B}_{\omega}^{\mathscr{F}}$ з індукованою напівгруповою операцією з  $\boldsymbol{B}_{\omega}^{\mathscr{F}_0}$ є піднапівгрупою в $\boldsymbol{B}_{\omega}^{\mathscr{F}_0}$.

У праці \cite{Gutik-Mykhalenych-2021} ми вивчаємо алгебраїчну структуру напівгрупи $\boldsymbol{B}_{\omega}^{\mathscr{F}}$ у випадку, коли ${\omega}$-замкнена сім'я $\mathscr{F}$ складається з індуктивних непорожніх підмножин у $\omega$, а саме групові конгруенції на напівгрупі $\boldsymbol{B}_{\omega}^{\mathscr{F}}$ та її гомоморфні ретракти. Доведено, що конгруенція $\mathfrak{C}$ на $\boldsymbol{B}_{\omega}^{\mathscr{F}}$ є груповою, тоді і лише тоді, коли звуження конгруенції $\mathfrak{C}$ на піднапівгрупу в $\boldsymbol{B}_{\omega}^{\mathscr{F}}$, яка ізоморфна біциклічній напівгрупі, не є відношенням рівності. Також, ми описуємо всі нетривіальні гомоморфні ретракти та ізоморфізми напівгрупи $\boldsymbol{B}_{\omega}^{\mathscr{F}}$.

Надалі скрізь в тексті ми вважаємо, що ${\omega}$-замкнена сім'я $\mathscr{F}$ складається лише з індуктивних непорожніх підмножин у $\omega$.

Очевидно, що кожний автоморфізм моноїда $\boldsymbol{B}_{\omega}^{\mathscr{F}}$ є тотожним перетворенням. У цій праці ми доводимо, що iн'\-cктивний ендоморфізм $\varepsilon$ напiвгрупи $\boldsymbol{B}_{\omega}^{\mathscr{F}}$ є тотожним відображенням тоді і тільки тоді, коли $\varepsilon$ має три різні нерухомі точки, що еквівалентно існуванню неідемпотентного елемента $(i,j,[p))\in\boldsymbol{B}_{\omega}^{\mathscr{F}}$ такого, що $(i,j,[p))\varepsilon=(i,j,[p))$.

%%%%%%%%%%%%%%%%%%%%%%%%%%%%%%%%%%%%%%%%%%%%%%%%%%%%%%%%%%%%%%%%%%%%%%%%%%%%%%%%%%
%\section{Про автоморфізми моноїда $\boldsymbol{B}_{\omega}^{\mathscr{F}}$}\label{section-2}

Надалі ми припускатимемо, що сім'я $\mathscr{F}$ містить щонайменше дві непорожні індуктивні підмножини у $\omega$. Також за тверд\-жен\-ням~1~\cite{Gutik-Mykhalenych-2021} для спрощення викладень не зменшуючи загальності можемо вважати, що $[0)\in\mathscr{F}$.

\begin{lemma}\label{lemma-2.1}
Якщо $\varepsilon$~--- ін'єктивний ендоморфізм моноїда $\boldsymbol{B}_{\omega}^{\mathscr{F}}$ такий, що $(1,1,[0))\varepsilon=(1,1,[0))$, то $\varepsilon$~--- тотожне відображення.
\end{lemma}

\begin{proof}
Доведемо твердження леми методом математичної індукції. Спочатку доведемо, що виконується початковий крок індукції: $(m,n,[0))\varepsilon=(m,n,[0))$ для довільних $m,n\in\omega$.
За тверд\-женням~3~\cite{Gutik-Mykhalenych-2020} для довільної множини $[n)\in \mathscr{F}$ піднапівгрупа $\boldsymbol{B}_{\omega}^{\{[n)\}}$ в $\boldsymbol{B}_{\omega}^{\mathscr{F}}$ ізоморфна біциклічній напівгрупі $\boldsymbol{B}_{\omega}$ стосовно відображення $(i,j,[n))\mapsto (i,j)$. Оскільки $(0,0,[0))$, $(1,1,[0))\in \boldsymbol{B}_{\omega}^{\{[0)\}}$, то за твердженням~4~\cite{Gutik-Mykhalenych-2021}, $(i,j,[0))\in \boldsymbol{B}_{\omega}^{\{[0)\}}$ справджується рівність $(i,j,[0))\varepsilon=(i,j,[0))$ для довільних $i,j\in\omega$. Позаяк $\boldsymbol{B}_{\omega}^{\mathscr{F}}$ --- інверсна напівгрупа та $(0,1,[0))^{-1}=(1,0,[0))$ в $\boldsymbol{B}_{\omega}^{\mathscr{F}}$, то з рівностей
\begin{align*}
  (0,1,[0))\cdot (1,0,[0))&=(0,0,[0)), \\
  (1,0,[0))\cdot (0,1,[0))&=(1,1,[0)), \\
  (0,0,[0))\varepsilon    &=(0,0,[0)), \\
  (1,1,[0))\varepsilon    &=(1,1,[0))
\end{align*}
і твердження~1.4.21$(1)$~\cite{Lawson-1998} випливає,  якщо $(0,1,[0))\varepsilon=(i,j,[0))$, то
\begin{equation*}
  (1,0,[0))\varepsilon=\big((0,1,[0))^{-1}\big)\varepsilon=\big((0,1,[0))\varepsilon\big)^{-1}=(i,j,[0))^{-1}=(j,i,[0)),
\end{equation*}
а отже,
\begin{equation*}
  (0,0,[0))=(0,0,[0))\varepsilon=\big((0,1,[0)){\cdot}(1,0,[0))\big)\varepsilon=(0,1,[0))\varepsilon{\cdot}(1,0,[0))\epsilon=(i,j,[0)){\cdot}(j,i,[0))= (i,i,[0))
\end{equation*}
i
\begin{equation*}
  (1,1,[0))=(1,1,[0))\varepsilon=\big((1,0,[0)){\cdot}(0,1,[0))\big)\varepsilon=(1,0,[0))\varepsilon{\cdot}(0,1,[0))\epsilon=(j,i,[0)){\cdot}(i,j,[0))= (j,j,[0)),
\end{equation*}
звідки впливає, що $(0,1,[0))\varepsilon=(0,1,[0))$ i $(1,0,[0))\varepsilon=(1,0,[0))$. Позаяк для довільних $m,n\in\omega$ в напівгрупі $\boldsymbol{B}_{\omega}^{\mathscr{F}}$ справджується рівність $(m,n,[0))=(1,0,[0))^m\cdot(0,1,[0))^n$, то з вище доведеного випливає, що $(m,n,[0))\varepsilon=(m,n,[0))$ для довільних $m,n\in\omega$.

Далі доведемо крок індукції: $(i,j,[k))\varepsilon=(i,j,[k))$ для довільних $i,j\in\omega$ i $[k)\in\mathscr{F}$.
Припустимо, що $[k+1)\in\mathscr{F}$ для деякого числа $k\in\omega$ і довільних $m,n\in\omega$ та невід'ємного цілого числа  $p\leqslant k$ виконується рівність $(m,n,[p))\varepsilon=(m,n,[p))$. Доведемо, що з цього припущення випливає рівність $(m,n,[k+1))\varepsilon=(m,n,[k+1))$ довільних $m,n\in\omega$.
З означення природного часткового порядку на напівґратці $E(\boldsymbol{B}_{\omega}^{\mathscr{F}})$ (див. твердження~2~\cite{Gutik-Mykhalenych-2021}) випливає, що $(0,0,[k+1))$ i $(1,1,[k-1))$~--- єдині такі ідемпотенти $e$ напівгрупи $\boldsymbol{B}_{\omega}^{\mathscr{F}}$, для яких виконується нерівність
\begin{equation*}
  (1,1,[k))\preccurlyeq e\preccurlyeq (0,0,[k)),
\end{equation*}
у випадку $k>0$,  а у випадку $k=0$ ідемпотент $(0,0,[k+1))$ єдиний, який задовольняє цю умову.
Тоді з твердження~1.4.21$(6)$~\cite{Lawson-1998} та припущення індукції випливає, що
\begin{equation*}
  (1,1,[k))=(1,1,[k))\varepsilon\preccurlyeq (0,0,[k+1))\varepsilon\preccurlyeq(0,0,[k))\varepsilon=(0,0,[k)),
\end{equation*}
а отже,
\begin{equation*}
  (0,0,[k+1))\varepsilon=(0,0,[k+1)).
\end{equation*}
Знову, оскільки $(1,1,[k+1))$ i $(2,2,[k-1))$~--- єдині такі ідемпотенти $e$ напівгрупи $\boldsymbol{B}_{\omega}^{\mathscr{F}}$, для яких виконується нерівність
\begin{equation*}
  (2,2,[k))\preccurlyeq e\preccurlyeq (1,1,[k)),
\end{equation*}
у випадку $k>0$,  а у випадку $k=0$ ідемпотент $(1,1,[k+1))$ єдиний, який задовольняє цю умову, то з твердження~1.4.21$(6)$~\cite{Lawson-1998} та припущення індукції  випливає, що
\begin{equation*}
  (2,2,[k))=(2,2,[k))\varepsilon\preccurlyeq (1,1,[k+1))\varepsilon\preccurlyeq(1,1,[k))\varepsilon=(1,1,[k)),
\end{equation*}
а отже,
\begin{equation*}
  (1,1,[k+1))\varepsilon=(1,1,[k+1)).
\end{equation*}
Далі, аналогічно як і у випадку $k=0$ доводиться, що $(m,n,[k+1))\varepsilon=(m,n,[k+1))$ для довільних $m,n\in\omega$.
\end{proof}

\begin{lemma}\label{lemma-2.2}
Нехай $\varepsilon$~--- ін'єктивний ендоморфізм моноїда $\boldsymbol{B}_{\omega}^{\mathscr{F}}$. Якщо $(0,0,[k))\varepsilon=(0,0,[k))$ для деякого натурального числа $k\geqslant2$, то $\varepsilon$~--- тотожне відображення.
\end{lemma}

\begin{proof}
Якщо $(0,0,[k))\varepsilon=(0,0,[k))$ для деякого натурального числа $k\geqslant2$, то з визначення природного часткового порядку на   $E(\boldsymbol{B}_{\omega}^{\mathscr{F}})$ (див. твердження~2~\cite{Gutik-Mykhalenych-2021}) випливає, що
\begin{equation*}
  {\uparrow_{\preccurlyeq}}(0,0,[k))=\left\{(0,0,[0)), (0,0,[1)),\ldots,(0,0,[k))\right\},
\end{equation*}
а отже, множина ${\uparrow_{\preccurlyeq}}(0,0,[k))$ є максимальним $(k+1)$-елементним ланцюгом, який містить ідемпотенти $(0,0,[0))$ i $(0,0,[k))$, як найбільший та найменший елементи, відповідно. Також, з визначення природного часткового порядку на напвіґратці  $E(\boldsymbol{B}_{\omega}^{\mathscr{F}})$ випливає, що для ідемпотента  $(j,j,[k))\in\boldsymbol{B}_{\omega}^{\mathscr{F}}$ множина ${\uparrow_{\preccurlyeq}}(j,j,[k))$ є ланцюгом в $E(\boldsymbol{B}_{\omega}^{\mathscr{F}})$ тоді і лише тоді, коли $j=0$. Отже, отримуємо, що $(0,0,[j))\varepsilon=(0,0,[j))$ для всіх $j=1,\ldots,k$.

Далі ми доведемо, що $(1,1,[0))\varepsilon=(1,1,[0))$. Оскільки $(0,0,[0))\varepsilon=(0,0,[0))$, то за твердженням~3 з \cite{Gutik-Mykhalenych-2020}  піднапівгрупа $\boldsymbol{B}_{\omega}^{\{[0)\}}$ в $\boldsymbol{B}_{\omega}^{\mathscr{F}}$ ізоморфна біциклічній напівгрупі $\boldsymbol{B}_{\omega}$, і тоді з умови $(0,0,[0))$, $(1,1,[0))\in \boldsymbol{B}_{\omega}^{\{[0)\}}$, твердження~4~\cite{Gutik-Mykhalenych-2021} та твердження~1.4.21$(2)$~\cite{Lawson-1998} випливає, що $(1,1,[0))\varepsilon$~--- ідемпотент напівгрупи $\boldsymbol{B}_{\omega}^{\{[0)\}}$. За лемою~2 з \cite{Gutik-Mykhalenych-2020}, $(1,1,[0))\varepsilon=(i,i,[0))$ для деякого натурального числа $i$, оскільки $\varepsilon$~--- ін'єктивний ендоморфізм.

Якщо $i\geqslant 2$, то з твердження~2~\cite{Gutik-Mykhalenych-2021} випливає, що
\begin{equation*}
(1,1,[0))\varepsilon=(i,i,[0))\preccurlyeq(0,0,[2))=(0,0,[2))\varepsilon,
\end{equation*}
тобто
\begin{equation}\label{eq-2.1}
(1,1,[0))\varepsilon\cdot(0,0,[2))\varepsilon=(1,1,[0))\varepsilon.
\end{equation}
Однак, ідемпотенти $(1,1,[0))$ i $(0,0,[2))$ за твердженням~2~\cite{Gutik-Mykhalenych-2021} непорівняльні в $E(\boldsymbol{B}_{\omega}^{\mathscr{F}})$, а отже
\begin{equation}\label{eq-2.2}
(0,0,[2))\neq(1,1,[0))\cdot(0,0,[2))\neq (1,1,[0)).
\end{equation}
Отримані умови \eqref{eq-2.1} i \eqref{eq-2.2} суперечать ін'єктивності ендоморфізму $\varepsilon$. З отриманого протиріччя випливає, що $i=1$, а отже, $(1,1,[0))\varepsilon=(1,1,[0))$. Далі скористаємося лемою~\ref{lemma-2.1}.
\end{proof}

З означення напівгрупової операції на $\boldsymbol{B}_{\omega}^{\mathscr{F}}$ випливає, що у випадку, коли $\mathscr{F}$ ---  ${\omega}$-замкнена сім'я підмножин у $\omega$ та $F\in\mathscr{F}$~--- непорожня індуктивна підмножина в $\omega$, то множина
\begin{equation*}
  \boldsymbol{B}_{\omega}^{\{F\}}=\left\{(i,j,F)\colon i,j\in\omega\right\}
\end{equation*}
з індукованою напівгруповою операцією з $\boldsymbol{B}_{\omega}^{\mathscr{F}}$ є піднапівгрупою в $\boldsymbol{B}_{\omega}^{\mathscr{F}}$, яка за тверд\-жен\-ням 3 з~\cite{Gutik-Mykhalenych-2020} ізоморфна біциклічній напівгрупі.

\begin{lemma}\label{lemma-2.3}
Нехай $\varepsilon$~--- ін'єктивний ендоморфізм моноїда $\boldsymbol{B}_{\omega}^{\mathscr{F}}$. Якщо існує такий ідемпотент \linebreak $(i,i,[p))\in\boldsymbol{B}_{\omega}^{\mathscr{F}}\setminus \{(0,0,[0)),(0,0,[1))\}$, що $(i,i,[p))\varepsilon=(i,i,[p))$, то $\varepsilon$~--- тотожне відображення.
\end{lemma}

\begin{proof}
Розглянемо можливі випадки:
\begin{itemize}
  \item[$(i)$] $i\geqslant 1$ i $p=0$;
  \item[$(ii)$] $i\geqslant 1$ i $p=1$;
  \item[$(iii)$] $p\geqslant 2$.
\end{itemize}

\smallskip

$(i)$ Якщо $i=1$ i $p=0$, то твердження леми випливає з леми~\ref{lemma-2.1}. Тому надалі будемо вважати, що $i\geqslant2$.

За твердженням~3 з \cite{Gutik-Mykhalenych-2020} піднапівгрупа $\boldsymbol{B}_{\omega}^{\{[0)\}}$ моноїда $\boldsymbol{B}_{\omega}^{\mathscr{F}}$ ізоморфна біциклічній напівгрупі $\boldsymbol{B}_{\omega}$, і оскільки $(i,i,[0))\varepsilon=(i,i,[0))$, то за наслідком~1.32 \cite{Clifford-Preston-1961} образ піднапівгрупи $\boldsymbol{B}_{\omega}^{\{[0)\}}$ стосовно ендоморфізму $\varepsilon$ ізоморфний біциклічній напівгрупі, а тоді  образ $(\boldsymbol{B}_{\omega}^{\{[0)\}})\varepsilon$  міститься в $\boldsymbol{B}_{\omega}^{\{[0)\}}$. З визначення природного часткового порядку на $\boldsymbol{B}_{\omega}^{\{[0)\}}$, рівності $(i,i,[0))\varepsilon=(i,i,[0))$ та ін'єктивності ендоморфізму $\varepsilon$ за твердженням~1.4.21$(6)$ з \cite{Lawson-1998} маємо, що $(j,j,[0))\varepsilon=(j,j,[0))$ для всіх $j=1,\ldots,i$. Далі скористаємося лемою~\ref{lemma-2.1}.

\smallskip

$(ii)$ Розглянемо випадок  $i\geqslant 1$ i $p=1$. За твердженням~3 з \cite{Gutik-Mykhalenych-2020} піднапівгрупа $\boldsymbol{B}_{\omega}^{\{[1)\}}$ в $\boldsymbol{B}_{\omega}^{\mathscr{F}}$ ізоморфна біциклічній напівгрупі $\boldsymbol{B}_{\omega}$, і оскільки $(i,i,[1))\varepsilon=(i,i,[1))$, то за наслідком~1.32 з \cite{Clifford-Preston-1961} образ піднапівгрупи $\boldsymbol{B}_{\omega}^{\{[1)\}}$ стосовно ін'єктивного ендоморфізму $\varepsilon$ ізоморфний біциклічній напівгрупі, а тоді цей образ  міститься в $\boldsymbol{B}_{\omega}^{\{[1)\}}$. Оскільки ${\uparrow_{\preccurlyeq}}(i,i,[1))\cap \boldsymbol{B}_{\omega}^{\{[1)\}}=\{(0,0,[1)),\ldots,(i,i,[1))\}$, то визначення природного часткового порядку на $\boldsymbol{B}_{\omega}^{\{[1)\}}$, рівності $(i,i,[1))\varepsilon=(i,i,[1))$ та ін'єктивності ендоморфізму $\varepsilon$ за твердженням~1.4.21$(6)$~\cite{Lawson-1998} маємо, що $(j,j,[1))\varepsilon=(j,j,[1))$ для всіх $j=0,\ldots,i$.

З визначення природного часткового порядку на $E(\boldsymbol{B}_{\omega}^{\mathscr{F}})$  випливає, що
\begin{equation*}
{\uparrow_{\preccurlyeq}}(1,1,[1))=\{(0,0,[0)),(1,1,[0)),(0,0,[1)),(1,1,[1))\}.
\end{equation*}
Оскільки за твердженням~1.4.21$(6)$~\cite{Lawson-1998} гомоморфізм інверсних напівгруп зберігає природний част\-ковий порядок, а отже, і на напівґратках їхніх ідемпотентів, то з рівностей
\begin{align*}
  (0,0,[0))\varepsilon&=(0,0,[0)), \\
  (0,0,[1))\varepsilon&=(0,0,[1)), \\
  (1,1,[1))\varepsilon&=(1,1,[1))
\end{align*}
отримуємо, що
\begin{equation*}
\begin{split}
  {\uparrow_{\preccurlyeq}}\big((1,1,[1))\varepsilon\big)
    &=\{(0,0,[0))\varepsilon,(1,1,[0))\varepsilon,(0,0,[1))\varepsilon,(1,1,[1))\varepsilon\}= \\
    &=\{(0,0,[0)),(1,1,[0))\varepsilon,(0,0,[1)),(1,1,[1))\},
\end{split}
\end{equation*}
звідки випливає, що $(1,1,[0))\varepsilon=(1,1,[0))$. Далі скористаємося лемою~\ref{lemma-2.1}.

\smallskip

$(iii)$ Якщо $i=0$ i $p\geqslant2$, то твердження леми випливає з леми~\ref{lemma-2.2}. Тому надалі будемо вважати, що $i\geqslant1$.

За твердженням~3 з \cite{Gutik-Mykhalenych-2020} піднапівгрупа $\boldsymbol{B}_{\omega}^{\{[p)\}}$ в $\boldsymbol{B}_{\omega}^{\mathscr{F}}$ ізоморфна біциклічній напівгрупі $\boldsymbol{B}_{\omega}$, і оскільки $(i,i,[p))\varepsilon=(i,i,[p))$, то за наслідком~1.32 \cite{Clifford-Preston-1961} образ піднапівгрупи $\boldsymbol{B}_{\omega}^{\{[p)\}}$ стосовно ін'єктивного ендоморфізму $\varepsilon$ ізоморфний біциклічній напівгрупі, а тоді цей образ  міститься в $\boldsymbol{B}_{\omega}^{\{[p)\}}$. Оскільки ${\uparrow_{\preccurlyeq}}(i,i,[p))\cap \boldsymbol{B}_{\omega}^{\{[p)\}}=\{(0,0,[p)),\ldots,(i,i,[p))\}$, то визначення природного часткового порядку на $\boldsymbol{B}_{\omega}^{\{[p)\}}$, рівності $(i,i,[p))\varepsilon=(i,i,[p))$ та ін'єктивності ендоморфізму $\varepsilon$ за твердженням~1.4.21$(6)$~\cite{Lawson-1998} маємо, що $(j,j,[p))\varepsilon=(j,j,[p))$ для всіх $j=0,\ldots,i$. Отже, $(0,0,[p))\varepsilon=(0,0,[p))$ і залишилося скористатися лемою~\ref{lemma-2.2}.
\end{proof}

\begin{theorem}\label{theorem-2.4}
Нехай $\mathscr{F}$~--- $\omega$-замкнена сім'я індуктивних непорожніх підмножин у $\omega$. Тоді ін'єк\-тив\-ний ендоморфізм $\varepsilon$ моноїда $\boldsymbol{B}_{\omega}^{\mathscr{F}}$ є тотожним відображенням тоді і тільки тоді, коли існує такий елемент $(i,j,[p))\in\boldsymbol{B}_{\omega}^{\mathscr{F}}\setminus \{(0,0,[0)),(0,0,[1))\}$, що $(i,j,[p))\varepsilon=(i,j,[p))$.
\end{theorem}

\begin{proof}
\textbf{Необхідність} очевидна.

Доведемо \textbf{достатність}. Якщо елемент  $(i,j,[p))\in\boldsymbol{B}_{\omega}^{\mathscr{F}}\setminus \{(0,0,[0)),(0,0,[1))\}$ є ідемпотентом, то твердження теореми випливає з леми~\ref{lemma-2.3}. Припустимо, що $(i,j,[p))$ не є ідемпотентом. Тоді з леми~2~\cite{Gutik-Mykhalenych-2020} випливає, що $i\neq j$, а з твердження~1.4.21$(1)$~\cite{Lawson-1998} отримуємо, що
\begin{equation*}
  (j,i,[p))\varepsilon=\big((i,j,[p))^{-1}\big)\varepsilon=\big((i,j,[p))\varepsilon\big)^{-1}=(i,j,[p))^{-1}=(j,i,[p)).
\end{equation*}
Отже, маємо, що
\begin{equation*}
  (j,j,[p))\varepsilon=\big((j,i,[p))\cdot(i,j,[p))\big)\varepsilon=(j,i,[p))\varepsilon\cdot(i,j,[p))\varepsilon=(j,i,[p))\cdot(i,j,[p))=(j,j,[p))
\end{equation*}
i
\begin{equation*}
  (i,i,[p))\varepsilon=\big((i,j,[p))\cdot(j,i,[p))\big)\varepsilon=(i,j,[p))\varepsilon\cdot(j,i,[p))\varepsilon=(i,j,[p))\cdot(j,i,[p))=(i,i,[p)),
\end{equation*}
і оскільки $i\neq j$, то виконується хоча б одна з умов $(i,i,[p))\in\boldsymbol{B}_{\omega}^{\mathscr{F}}\setminus \{(0,0,[0)),(0,0,[1))\}$, або $(j,j,[p))\in\boldsymbol{B}_{\omega}^{\mathscr{F}}\setminus \{(0,0,[0)),(0,0,[1))\}$. Залишилося скористатися лемою~\ref{lemma-2.3}.
\end{proof}

%З доведення теореми~\ref{theorem-2.4} випливає така теорема.

\begin{theorem}\label{theorem-2.5}
Нехай $\mathscr{F}$~--- $\omega$-замкнена сім'я індуктивних непорожніх підмножин у $\omega$, яка містить хоча б дві множини. Тоді для ін'єктивного ендоморфізму $\varepsilon$ моноїда $\boldsymbol{B}_{\omega}^{\mathscr{F}}$ такі умови еквівалентні:
\begin{itemize}
  \item[$(i)$]   $\varepsilon$~--- тотожне відображення;
  \item[$(ii)$]  існує неідемпотентний елемент $(i,j,[p))\in\boldsymbol{B}_{\omega}^{\mathscr{F}}$ такий, що $(i,j,[p))\varepsilon=(i,j,[p))$.
  \item[$(iii)$] $\varepsilon$ має хоча б три нерухомі елементи.
\end{itemize}
\end{theorem}

\begin{proof}
Імплікації $(i)\Rightarrow(ii)$ та $(i)\Rightarrow(iii)$ очевидні. Імплікація $(ii)\Rightarrow(i)$ з теореми~\ref{theorem-2.4}.

$(iii)\Rightarrow(i)$ Очевидно, що $(0,0,[0))\varepsilon=(0,0,[0))$. Якщо один з нерухомих елементів елементів $\varepsilon$ неідемпотентний, то виконується умова $(ii)$. Тому припустимо, що $(i,i,[p))\varepsilon=(i,i,[p))$ та $(j,j,[q))\varepsilon=(j,j,[q))$ для двох різних неодиничних ідемпотентів $(i,i,[p))$ і $(j,j,[q))$ моноїда $\boldsymbol{B}_{\omega}^{\mathscr{F}}$.

Тоді один з ідемпотентів $(i,i,[p))$ або $(j,j,[q))$ відмінний від $(0,0,1)$. Далі скористаємося теоремою~\ref{theorem-2.4}.
\end{proof}

\begin{remark}\label{remark-2.6}
За твердженням~3 з~\cite{Gutik-Mykhalenych-2020} для $\omega$-замкненої сім'ї $\mathscr{F}$ індуктивних непорожніх підмножин у $\omega$ напівгрупа $\boldsymbol{B}_{\omega}^{\mathscr{F}}$ ізоморфна біциклічному моноїду $\boldsymbol{B}_{\omega}$ тоді і лише тоді, коли сім'я $\mathscr{F}$ скла\-да\-єть\-ся з єдиної непорожньої індуктивної підмножини в $\omega$. Усі автоморфізми біциклічного моноїда тривіальні, а напівгрупа $\mathrm{\mathbf{End}}(\boldsymbol{B}_{\omega})$  ендоморфізмів  бі\-цик\-ліч\-ної напівгрупи $\boldsymbol{B}_{\omega}$ ізоморфна напівпрямому добутку $(\omega,+)\rtimes_\varphi(\omega,*)$, де $(\omega,+)$ і $(\omega,*)$ --- адитивна та мультиплікативна напівгрупи невід'ємних цілих чисел, відповідно (див. \cite{Gutik-Prokhorenkova-Sekh-2021}). Оскільки біциклічний моноїд $\boldsymbol{B}_{\omega}$ як інверсна напівгрупа породжується єдиним елементом $(0,1)$, то очевидно, що для його ін'єктивного  ендоморфізму $\varepsilon$  такі умови еквівалентні:
\begin{itemize}
  \item[$(i)$]   $\varepsilon$~--- тотожне відображення;
  \item[$(ii)$]  існує неідемпотентний елемент $(i,j)\in\boldsymbol{B}_{\omega}$ такий, що $(i,j)\varepsilon=(i,j)$.
  \item[$(iii)$] $\varepsilon$ має хоча б два нерухомі елементи.
\end{itemize}
\end{remark}

У випадку ендоморфізму біциклічної напівгрупи $\boldsymbol{B}_{\omega}$ умову $(ii)$ в зауваженні \ref{remark-2.6}  можна замінити твердженням~\ref{proposition-2.6}.

\begin{proposition}\label{proposition-2.6}
Нехай $\varepsilon$~--- ін'єктивний ендоморфізм біцикліної напівгрупи $\boldsymbol{B}_{\omega}$. Якщо $(i,j)\varepsilon=(i,j)$ для деякого неодиничного елемента $(i,j)\in \boldsymbol{B}_{\omega}$, то $\varepsilon$~--- тотожне перетворення.
\end{proposition}

\begin{proof}
Спочатку розглянемо випадок, коли $(i,i)\varepsilon=(i,i)$ для деякого неодиничного ідемпотента $(i,i)\in \boldsymbol{B}_{\omega}$. З ін'єктивності ендоморфізму $\varepsilon$ та з того, що гомоморфізм інверсних напівгруп зберігає природний частковий порядок і образ ідемпотента при гомоморфізмі напівгруп є знову ідемпотент (див. твердження~1.4.21 з~\cite{Lawson-1998}) випливає, що
$(k,k)\varepsilon=(k,k)$ для всіх $k=0,\ldots,i-1$, оскільки ${\uparrow_{\preccurlyeq}}(i,i)=\{(k,k)\colon k=0,\ldots,i\}$. Далі скористаємося еквівалентністю умов $(i)$ та $(iii)$ зауваження~\ref{remark-2.6}.

Далі припустимо, що $(i,j)\in \boldsymbol{B}_{\omega}$ не є ідемпотентом.
За твердженням~1.4.21 з~\cite{Lawson-1998} маємо, що
\begin{align*}
  (i,i)\varepsilon&=((i,j)\cdot(j,i))\varepsilon=((i,j)\cdot(i,j)^{-1})\varepsilon= (i,j)\varepsilon\cdot((i,j)^{-1})\varepsilon= (i,j)\varepsilon\cdot((i,j)\varepsilon)^{-1}=\\
   &=(i,j)\cdot(i,j)^{-1}=(i,j)\cdot(j,i)=(i,i),
\end{align*}
та аналогічно, $(j,j)\varepsilon=(j,j)$.

Оскільки $(i,j)\neq(0,0)$, то $(i,i)\neq(j,j)$. Отже, $(i,i)\varepsilon=(i,i)$ для деякого неодиничного ідемпотента $(i,i)\in \boldsymbol{B}_{\omega}$. Далі скористаємося попередньою частиною доведення.
\end{proof}

Тому природно виникає таке запитання.

\begin{question}\label{question-2.7}
Чи можна послабити умови теореми~\ref{theorem-2.5}?
\end{question}

З прикладу \ref{example-2.8} випливає, що умову про ін'єктивність ендоморфізму $\varepsilon$ моноїда $\boldsymbol{B}_{\omega}^{\mathscr{F}}$ не можна вилучити в припущеннях теореми~\ref{theorem-2.5}.

\begin{example}\label{example-2.8}
Нехай $\mathscr{F}$~--- довільна $\omega$-замкнена сім'я індуктивних непорожніх підмножин у $\omega$, яка містить хоча б дві множини. Не зменшуючи загальності можемо вважати, що $[0),[1)\in\mathscr{F}$. Означимо відображення $\mu_0\colon \boldsymbol{B}_{\omega}^{\mathscr{F}}\to \boldsymbol{B}_{\omega}^{\mathscr{F}}$ за формулою
\begin{equation*}
 (i,j,[p))\mu_0=(i,j,[0)), \qquad \hbox{для довільних} \quad i,j\in\omega \quad \hbox{i} \quad [p)\in\mathscr{F}.
\end{equation*}
За твердженням~5~\cite{Gutik-Mykhalenych-2021} відображення $\mu_0$ є ендоморфізмом напівгрупи $\boldsymbol{B}_{\omega}^{\mathscr{F}}$, причому піднапівгрупа $\boldsymbol{B}_{\omega}^{\{[0)\}}$ є гомоморфним ретрактом моноїда $\boldsymbol{B}_{\omega}^{\mathscr{F}}$ відносно ендоморфізму $\mu_0$, а отже, виконуються умови $(ii)$ i $(iii)$ теореми~\ref{theorem-2.5}, і не виконується умова $(i)$.
\end{example}

З прикладу \ref{example-2.9} випливає, що умову про те, що ендоморфізм $\varepsilon$ моноїда $\boldsymbol{B}_{\omega}^{\mathscr{F}}$ має хоча б три нерухомі елементи не можна послабити в умові $(iii)$ теореми~\ref{theorem-2.5}.

\begin{example}\label{example-2.9}
Нехай $\mathscr{F}=\{[0),[1)\}$. Зафіксуємо довільне натуральне число $k\geqslant 2$. Означимо відоб\-раження $\mu_k\colon \boldsymbol{B}_{\omega}^{\mathscr{F}}\to \boldsymbol{B}_{\omega}^{\mathscr{F}}$ за формулою
\begin{equation*}
  (i,j,[p))\mu_k=(ki,kj,[p)).
\end{equation*}

Очевидно, що так означене відображення $\mu_k$ ін'єктивне. Доведемо, що $\mu_k\colon \boldsymbol{B}_{\omega}^{\mathscr{F}}\to \boldsymbol{B}_{\omega}^{\mathscr{F}}$ --- гомоморфізм. Надалі будемо вважати, що $p\in\{0.1\}$. Тоді маємо, що
\begin{align*}
  \big((i_1,j_1,[p))\cdot(i_2,j_2,[0))\big)\mu_k&=
    \left\{
    \begin{array}{ll}
      (i_1-j_1+i_2,j_2,(j_1-i_2+[p))\cap [0))\mu_k, & \hbox{якщо~} j_1<i_2;\\
      (i_1,j_2,[p)\cap [0))\mu_k,                   & \hbox{якщо~} j_1=i_2;\\
      (i_1,j_1-i_2+j_2,[p)\cap (i_2-j_1+[0)))\mu_k, & \hbox{якщо~} j_1>i_2
    \end{array}
  \right.=
     \\
   &=
   \left\{
    \begin{array}{ll}
      (i_1-j_1+i_2,j_2,[0))\mu_k, & \hbox{якщо~} j_1<i_2;\\
      (i_1,j_2,[p))\mu_k,         & \hbox{якщо~} j_1=i_2;\\
      (i_1,j_1-i_2+j_2,[p))\mu_k, & \hbox{якщо~} j_1>i_2
    \end{array}
  \right.=
     \\
   &=
   \left\{
    \begin{array}{ll}
      (k(i_1-j_1+i_2),kj_2,[0)), & \hbox{якщо~} j_1<i_2;\\
      (ki_1,kj_2,[p)),           & \hbox{якщо~} j_1=i_2;\\
      (ki_1,k(j_1-i_2+j_2),[p)), & \hbox{якщо~} j_1>i_2
    \end{array}
  \right.
\end{align*}
\begin{align*}
  (i_1,j_1,[p))\mu_k\cdot(i_2,j_2,[0))\mu_k&=(ki_1,kj_1,[p))\cdot(ki_2,kj_2,[0))= \\
   &=
   \left\{
    \begin{array}{ll}
      (ki_1-kj_1+ki_2,kj_2,(kj_1-ki_2+[p))\cap [0)), & \hbox{якщо~} kj_1<ki_2;\\
      (ki_1,kj_2,[p)\cap [0)),                       & \hbox{якщо~} kj_1=ki_2;\\
      (ki_1,kj_1-ki_2+kj_2,[p)\cap (ki_2-kj_1+[0))), & \hbox{якщо~} kj_1>ki_2
    \end{array}
  \right.=
  \\
   &=
   \left\{
    \begin{array}{ll}
      (k(i_1-j_1+i_2),kj_2,[0)), & \hbox{якщо~} j_1<i_2;\\
      (ki_1,kj_2,[p)),           & \hbox{якщо~} j_1=i_2;\\
      (ki_1,k(j_1-i_2+j_2),[p)), & \hbox{якщо~} j_1>i_2
    \end{array}
  \right.
\end{align*}
i
\begin{align*}
  \big((i_1,j_1,[p))\cdot(i_2,j_2,[1))\big)\mu_k&=
    \left\{
    \begin{array}{ll}
      (i_1-j_1+i_2,j_2,(j_1-i_2+[p))\cap [1))\mu_k, & \hbox{якщо~} j_1<i_2;\\
      (i_1,j_2,[p)\cap [1))\mu_k,                   & \hbox{якщо~} j_1=i_2;\\
      (i_1,j_1-i_2+j_2,[p)\cap (i_2-j_1+[1)))\mu_k, & \hbox{якщо~} j_1>i_2
    \end{array}
  \right.=
     \\
   &=
   \left\{
    \begin{array}{ll}
      (i_1-j_1+i_2,j_2,[1))\mu_k, & \hbox{якщо~} j_1<i_2;\\
      (i_1,j_2,[1))\mu_k,         & \hbox{якщо~} j_1=i_2;\\
      (i_1,j_1-i_2+j_2,[p))\mu_k, & \hbox{якщо~} j_1>i_2
    \end{array}
  \right.=
     \\
   &=
   \left\{
    \begin{array}{ll}
      (k(i_1-j_1+i_2),kj_2,[1)), & \hbox{якщо~} j_1<i_2;\\
      (ki_1,kj_2,[1)),           & \hbox{якщо~} j_1=i_2;\\
      (ki_1,k(j_1-i_2+j_2),[p)), & \hbox{якщо~} j_1>i_2,
    \end{array}
  \right.
\end{align*}
\begin{align*}
  (i_1,j_1,[p))\mu_k\cdot(i_2,j_2,[1))\mu_k&=(ki_1,kj_1,[p))\cdot(ki_2,kj_2,[1))= \\
   &=
   \left\{
    \begin{array}{ll}
      (ki_1-kj_1+ki_2,kj_2,(kj_1-ki_2+[p))\cap [1)), & \hbox{якщо~} kj_1<ki_2;\\
      (ki_1,kj_2,[p)\cap [1)),                       & \hbox{якщо~} kj_1=ki_2;\\
      (ki_1,kj_1-ki_2+kj_2,[p)\cap (ki_2-kj_1+[1))), & \hbox{якщо~} kj_1>ki_2
    \end{array}
  \right.=
  \\
   &=
   \left\{
    \begin{array}{ll}
      (k(i_1-j_1+i_2),kj_2,[1)), & \hbox{якщо~} j_1<i_2;\\
      (ki_1,kj_2,[1)),           & \hbox{якщо~} j_1=i_2;\\
      (ki_1,k(j_1-i_2+j_2),[p)), & \hbox{якщо~} j_1>i_2,
    \end{array}
  \right.
\end{align*}
а отже, $\mu_k\colon \boldsymbol{B}_{\omega}^{\mathscr{F}}\to \boldsymbol{B}_{\omega}^{\mathscr{F}}$ --- ін'єктивний ендоморфізм.
\end{example}

%%%%%%%%%%%%%%%%%%%%%%%%%%%%%%%%%%%%%%%%%%%%%%%%%%%%%%%%%%%%%%%%%%%%%%%%%%%%%%%%%%
%\bigskip

%\section*{\textbf{Подяка}}

%Автори висловлюють щиру подяку  рецензентові за цінні поради та зауваження.

%%%%%%%%%%%%%%%%%%%%%%%%%%%%%%%%%%%%%%%%%%%%%%%%%%%%%%%%%%%%%%%%%%%%%%%%%%%%

%\vskip1cm

\begin{thebibliography}{10}

\bibitem{Vagner-1952}
В. В.~Вагнер,
\emph{Обощенные группы},
ДАН СССР \textbf{84} (1952), 1119--1122.

\bibitem{Gutik-Mykhalenych-2020}
О.~Гутік, М. Михаленич,
\emph{Про одне узагальнення бiциклiчного моноїда},
Вісник Львів. ун-ту. Сер. мех.-мат. \textbf{90} (2020), 5--19.

\bibitem{Gutik-Mykhalenych-2021}
О.~Гутік, М. Михаленич,
\emph{Про груповi конгруенцiї на напiвгрупi $\boldsymbol{B}_{\omega}^{\mathscr{F}}$ та \"{\i}\"{\i} гомоморфнi ретракти у випадку, коли сiм'я $\mathscr{F}$ складається з непорожнiх iндуктивних пiдмножин у $\omega$},
Вісник Львів. ун-ту. Сер. мех.-мат. \textbf{91} (2021), 5--27.


\bibitem{Gutik-Pozdnyakova-2021??}
O. В. Гутік, І. В. Позднякова,
\emph{Про напiвгрупу, породжену розширеною бiциклiчною напiвгрупою та $\omega$-замкненою сiм'єю},
Мат. методи  фіз.-мех. поля \textbf{64} (2021), no.~1, 21--34 (arXiv:2107.14075).


\bibitem{Gutik-Prokhorenkova-Sekh-2021}
О. Гутік, О. Прохоренкова, Д. Сех,
\emph{Про ендоморфiзми бiциклiч\-ної напiвгрупи та розширеної бiциклiчної напiвгрупи},
Вісник Львів. ун-ту. Сер. мех.-мат. \textbf{92} (2021), 5--16.

%\bibitem{Gutik-Savchuk-2018}
%О. Гутік, А. Савчук,
%\emph{Напiвгрупа часткових коскiнченних iзометрiй натуральних чисел},
%Буковинський мат. журнал \textbf{6} (2018), no. 1--2, 42--51.

%\bibitem{Andersen-1952}
%O. Andersen,
%\emph{Ein Bericht uber die Struk\-tur abstrakter Halbgruppen},
%PhD Thesis. Ham\-burg, 1952.

%\bibitem{Ash-1979}
%C. J. Ash,
%\emph{The $\mathscr{J}$-classes of an inverse semigroup},
%J. Austral. Math. Soc. Ser. A \textbf{28} (1979), 427--432.

%\bibitem{Brown-1965}
% D. R. Brown,
%\emph{Topological semilattices on the two-cell},
%Pacific J. Math. 15 (1965), no. 1, 35--46.


%\bibitem{Bruck-1958}
%R. H. Bruck,
%\emph{A survey of binary systems},
%(Erg. Math. Grenzgebiete. Neue Folge. Heft 20) Springer, Berlin-G\"{o}ttingen-Heidelberg,  1958.

\bibitem{Clifford-Preston-1961}
A.~H.~Clifford and G.~B.~Preston,
\emph{The Algebraic Theory of Semigroups}, Vol. I,
Amer. Math. Soc. Surveys {\bf 7}, Pro\-vi\-den\-ce, R.I.,  1961.


\bibitem{Clifford-Preston-1967}
A.~H.~Clifford and G.~B.~Preston,
\emph{The Algebraic Theory of Semigroups}, Vol.  II,
Amer. Math. Soc. Surveys {\bf 7}, Pro\-vi\-den\-ce, R.I.,   1967.

%\bibitem{Fortunatov-1976}
%V. A. Fortunatov,
%\emph{Congruences on simple extensions of semigroups},
%Semigroup Forum \textbf{13} (1976), 283--295.

%\bibitem{Fotedar-1974}
%G.~L.~Fotedar, \emph{On a semigroup associated with an ordered
%group}, Math. Nachr. \textbf{60} (1974), 297--302.

%\bibitem{Fotedar-1978}
%G.~L.~Fotedar, \emph{On a class of bisimple inverse semigroups},
%Riv. Mat. Univ. Parma (4) \textbf{4} (1978), 49--53.

%\bibitem{Green-1951}
%J. A. Green
%\emph{On the structure of semigroups},
%Ann. Math. Ser. 2 \textbf{54} (1951), no.~1, 163--172.

%\bibitem{Gutik-2018}
%O. Gutik,
%\emph{On locally compact semitopological $0$-bisimple inverse $\omega$-semigroups,}
%Topol. Algebra Appl. \textbf{6} (2018), 77--101.

\bibitem{Gutik-Lysetska=2021}
O. Gutik and O. Lysetska,
\emph{On the semigroup $\boldsymbol{B}_{\omega}^{\mathscr{F}}$ which is generated by the family $\mathscr{F}$ of atomic subsets of $\omega$},
Вісник Львів. ун-ту. Сер. мех.-мат. \textbf{92} (2021), 34--50 (arXiv:2108.11354).

%\bibitem{Gutik-Pagon-Pavlyk=2011}
%O. Gutik, D. Pagon, and K. Pavlyk,
%\emph{Congruences on bicyclic extensions of a linearly ordered group},
%Acta Comment. Univ. Tartu. Math. \textbf{15} (2011), no. 2, 61--80.

\bibitem{Gutik-Pozdniakova=2023}
O. Gutik and I. Pozdniakova,
\emph{On the group of automorphisms of the semigroup $\boldsymbol{B}^{\mathscr{F}}_{\mathbb{Z}}$ with the family $\mathscr{F}$ of inductive nonempty subsets of $\omega$},
Algebra Discrete Math. \textbf{35} (2023), no.~1 (to appear) (arXiv:2206.12819).

\bibitem{Lawson-1998}
M.~Lawson,
\emph{Inverse Semigroups. The Theory of Partial Symmetries},
World Scientific, Singapore, 1998.


\bibitem{Lysetska=2020}
O. Lysetska,
\emph{On feebly compact topologies on the semigroup $\boldsymbol{B}_{\omega}^{\mathscr{F}_1}$},
Вісник Львів. ун-ту. Сер. мех.-мат. \textbf{90} (2020), 48--56.



\bibitem{Petrich-1984}
M.~Petrich,
\emph{Inverse Semigroups},
John Wiley $\&$ Sons, New York, 1984.



\end{thebibliography}
\end{document}